\documentstyle[amssymb,12pt]{article}

\newtheorem{th}{Theorem}[section]
\newtheorem{lem}[th]{Lemma}
\newtheorem{prop}[th]{Proposition}
\newtheorem{cor}[th]{Corollary}
\newtheorem{defn}[th]{Definition}
\newenvironment{defn-new}{\begin{defn} \em}{\end{defn}}
\newtheorem{rem}[th]{Remark}
\newenvironment{rem-new}{\begin{rem} \em}{\end{rem}}
\newtheorem{ex}[th]{Example}
\newenvironment{ex-new}{\begin{ex} \em}{\end{ex}}

\newenvironment{notation-new}{\begin{rem} \em}{\end{rem}}

\newenvironment{agr-new}{\begin{rem} \em}{\end{rem}}

\makeatletter \@addtoreset{equation}{section} \makeatother

\makeatletter \@addtoreset{figure}{section} \makeatother

\begin{document}

\begin{center}
{\Large Einstein doubly warped product manifolds with a semi-symmetric
metric connection}

\bigskip {\large {\bf Punam Gupta and Abdoul Salam Diallo }}
\end{center}

\noindent {\bf {Abstract.} }In this paper, we study the doubly warped
product manifolds{\bf \ }with semi-symmetric metric connection. We derive
the curvatures formulas for doubly warped product manifold with
semi-symmetric metric connection in terms of curvatures of components of
doubly warped product manifolds. We also prove the necessary and sufficient
condition for a doubly warped product manifold to be a warped product
manifold. We obtain some results for Einstein doubly warped product manifold
and Einstein-like doubly warped product manifold of class ${\cal A}$ with
respect to a semi-symmetric metric connection.

\noindent {\bf 2010 Mathematics Subject Classification.} 53C05, 53C25,
53C50.\medskip

\noindent {\bf {Keywords.} }Doubly warped products, semi-symmetric metric
connection, Levi-Civita connection, Einstein manifolds, Einstein-like
manifold of class ${\cal A}$.

\section{Introduction\label{sect-intro}}

In 1969, Bishop and O' Neill \cite{Bishop-Neill-69} introduced singly warped
products or warped products. They used this concept to construct Riemannian
manifolds with negative sectional curvature. The warped product $B\times
_{h}F$ of two Riemannian manifolds $(B,g_{B})$ and $(F,g_{F})$ with a smooth
function $h:B\rightarrow (0,\infty )$ is a product manifold of form $B\times
F$ with the metric tensor $g=g_{B}\oplus h^{2}g_{F}$. Here, $(B,g_{B})$ is
called the base manifold, $(F,g_{F})$ is called as the fiber manifold and $h$
is called as the warping function. In 1983, O'Neill \cite{ONeill-83}
discussed warped products and derived curvature formulas of warped products
in terms of curvatures of components of warped products. After that, theory
of relativity demands a larger class of manifold and then idea of doubly
warped products was introduced. The doubly warped product $B_{f}\times
{}_{h}F$ of two Riemannian manifolds $(B,g_{B})$ and $(F,g_{F})$ with smooth
functions, which are known as warping functions, $h:B\rightarrow (0,\infty )$%
, $f:F\rightarrow (0,\infty )$ is a product manifold of form $B\times F$
with the metric tensor $g=f^{2}g_{B}\oplus h^{2}g_{F}$ . Doubly warped
products are studied by many authors: Allison \cite{Allison-91} studied
pseudoconvexity of Lorentzian doubly warped products; Gebarowski considered
conformal properties of doubly warped products in (\cite{Gebarowski-94},
\cite{Gebarowski-96}). Geodesic completeness of Riemannian and Lorentzian
doubly warped products are studied by Unal \cite{Unal-01}. In 2003, Ramos et
al. \cite{Ramos-03} gave an invariant characterization of doubly warped
space-times in terms of Newman-Penrose formalism, and proposed a
classification scheme. Recently, the first author \cite{Gupta} studied the
compact Einstein doubly warped product manifold and find out some
interesting results.

In 1924, Friedmann and Schouten \cite{Friedmann-Schouten-24} introduced the
idea of a semi-symmetric linear connection on a Riemannian manifold. Later,
in 1932, Hayden \cite{Hayden-32} introduced the concept of semi-symmetric
metric connection. A metric connection is Levi-Civita connection when its
torsion is zero and it becomes the Hayden connection \cite{Hayden-32} when
it has non-zero torsion. Thus, metric connections include both the
Levi-Civita connections and the Hayden connections. In 1970, Yano \cite%
{Yano-70} considered semi-symmetric metric connection and studied some of
its properties.

In 2001, Sular and \"{O}zgur \cite{Ozgur-Sular} studied the warped product
with semi-symmetric metric connection and find some results. Motivated by
these circumstances, we study doubly warped product manifolds with a
semi-symmetric metric connection and derive the curvatures formulas for
doubly warped product manifolds with semi-symmetric metric connection in
terms of curvatures of components of doubly warped product manifolds. We
also find results for Einstein doubly warped product manifold and
Einstein-like doubly warped product manifold of class ${\cal A}$ with
respect to a semi-symmetric metric connection.

\section{Preliminaries\label{sect-prelim}}

Let $M$ be an $n$-dimensional Riemannian manifold with Riemannian metric $g$%
. A linear connection $\tilde{\nabla}$ on a Riemannian manifold $M$ is
called a {\em semi-symmetric connection} if the torsion tensor $\tilde{T}$
of the connection $\tilde{\nabla}$ given by
\begin{equation}
\tilde{T}(X,Y)=\tilde{\nabla}_{X}Y-\tilde{\nabla}_{Y}X-[X,Y]  \label{torsion}
\end{equation}%
satisfies%
\begin{equation}
\tilde{T}(X,Y)=\pi (Y)X-\pi (X)Y,  \label{torssmc}
\end{equation}%
where $\pi $ is a $1$-form associated with the vector field $P$ on $M$
defined by%
\begin{equation}
\pi (X)=g(X,P).  \label{pi}
\end{equation}%
A semi-symmetric connection $\tilde{\nabla}$ is called a {\em semi-symmetric
metric connection} if $\tilde{\nabla}g=0$. Let $\nabla $ be the Levi-Civita
connection of a Riemannian manifold, then the unique semi-symmetric metric
connection $\tilde{\nabla}$ given by Yano \cite{Yano-70} is
\begin{equation}
\tilde{\nabla}_{X}Y=\nabla _{X}Y+\pi (Y)X-g(X,Y)P.  \label{semsymmetric}
\end{equation}%
\medskip

A relation between the curvature tensors $R$ and $\tilde{R}$ of the
Levi-Civita connection $\nabla $ and the semi-symmetric connection $\tilde{%
\nabla}$ is given by
\begin{eqnarray}
\tilde{R}(X,Y)Z &=&R(X,Y)Z+g(Z,\nabla _{X}P)Y-g(Z,\nabla _{Y}P)X  \nonumber
\\
&&+g(X,Z)\nabla _{Y}P-g(Y,Z)\nabla _{X}P  \nonumber \\
&&+\pi (P)\left( g(X,Z)Y-g(Y,Z)X\right)  \nonumber \\
&&+\left( g(Y,Z)\pi (X)-g(X,Z)\pi (Y)\right) P  \nonumber \\
&&+\pi (Z)\left( \pi (Y)X-\pi (X)Y\right) ,\quad X,Y,Z\in {\frak X}(M),
\label{curvature}
\end{eqnarray}%
where ${\frak X}(M)$ is the set of smooth vector fields on $M$ \cite{Yano-70}%
. Yano proved that a Riemannian manifold is conformally flat if and only if
it admits a semi-symmetric metric connection whose curvature tensor vanishes
identically. This result was generalized for vanishing Ricci tensor of the
semi-symmetric metric connection by T. Imai (see \cite{Imai-23}, \cite%
{Imai-23-Tensor}). For a general survey of different kinds of connections
(see Tripathi \cite{Tripathi08}).

\section{Doubly Warped Product Manifolds\label{warproduct}}

\noindent Doubly warped product manifolds were introduced as a
generalization of warped product manifolds. Let $(B,g_{B})$ and $(F,g_{F})$
be two Riemannian manifolds with real dimension $n_{1}$ and $n_{2}$,
respectively. Let $h:B\rightarrow
\mathbb{R}
^{+}$ and $f:F\rightarrow
\mathbb{R}
^{+}$ be two smooth functions. Consider the product manifold $B\times F$
with its projections $\rho :B\times F\rightarrow B$ and $\sigma :B\times
F\rightarrow F$. The {\em doubly warped product} $B_{f}\times {}_{h}F$ is
the product manifold $B\times F$ furnished with the metric tensor
\[
g=(f\circ \sigma )^{2}\rho ^{\ast }\left( g_{B}\right) +(h\circ \rho
)^{2}\sigma ^{\ast }\left( g_{F}\right) ,
\]%
where $^{\ast }$ denotes pullback. If $X$ is tangent to $B\times F$ at $%
(p,q) $, then
\[
g(X,X)=f^{2}(q)g_{B}(d\rho (X),d\rho (X))+h^{2}(p)g_{F}(d\sigma (X),d\sigma
(X)).
\]%
Thus, we have%
\begin{equation}
g=f^{2}g_{B}+h^{2}g_{F}.  \label{metric}
\end{equation}%
The functions $h$ and $f$ are called the warping functions of the doubly
warped product. The manifold $B$ is known as the base of $(M,g)$ and the
manifold $F$ is known as the fibre of $(M,g)$. If the warping function $f$
or $h$ is constant, then the doubly warped product $B_{f}\times {}_{h}F$
reduces to a warped product $B\times {}_{h}\tilde{F}$ (or $\tilde{B}%
_{f}\times F)$, where the fibre $\tilde{F}$ is just $F$ with metric $\tilde{g%
}_{F}$ given by $\frac{1}{h^{2}}g_{F}$ (or where the base $\tilde{B}$ is
just $B$ with metric $\tilde{g}_{B}$ given by $\frac{1}{f^{2}}g_{B}$). If
the warping function $f$ or $h$ is equal to $1$, then the doubly warped
product $B_{f}\times {}_{h}F$ reduces to a warped product $B\times {}_{h}F$
(or $B_{f}\times F$). If both $f$ and $h$ are constant, then it is simply a
product manifold $\tilde{B}\times \tilde{F}$, where the fibre $\tilde{F}$ is
just $F$ with metric $\tilde{g}_{F}$ given by $\frac{1}{h^{2}}g_{F}$ and the
base $\tilde{B}$ is just $B$ with metric $\tilde{g}_{B}$ given by $\frac{1}{%
f^{2}}g_{B}$. If both $f$ and $h$ are equal to $1$, then it is simply a
product manifold $B\times F$.

The set of all smooth and positive valued functions $h:B\rightarrow
\mathbb{R}
^{+}$ and $f:F\rightarrow
\mathbb{R}
^{+}$are denoted by ${\cal F}(B)=C^{\infty }(B)$ and ${\cal F}(F)=C^{\infty
}(F)$, respectively. The lift of $h$ and $f$ to $M$ are defined by $\tilde{h}%
=h\circ \rho \in $ ${\cal F}(M)$ and $\tilde{f}=f\circ \sigma \in $ ${\cal F}%
(M)$, respectively. If $X_{p}\in T_{p}(B)$ and $q\in F$, then the lift $%
\tilde{X}_{(p,q)}$ of $X_{p}$ to $M$ is the unique tangent vector in $%
T_{(p,q)}(B\times \{q\})$ such that $d\rho _{(p,q)}(\tilde{X}_{(p,q)})=X_{p}$
and $d\sigma _{(p,q)}(\tilde{X}_{(p,q)})=0$. The set of all such horizontal
tangent vector lifts will be denoted by $L_{(p,q)}(B)$. Similarly, we can
define the set of all vertical tangent vector lifts $L_{(p,q)}(F)$.

Let $X\in {\frak X}(B)$, where ${\frak X}(B)$ is the set of smooth vector
fields on $B$. The lift $\tilde{X}$ of $X$ to $M$ is the unique element of $%
{\frak X}(M)$ whose value at each $(p,q)$ is the lift of $X_{p}$ to $(p,q)$.
The set of such lifts will be denoted by ${\cal L}(B)$. In a similar manner,
we can define ${\cal L}(F)$. Throughout the paper, we assume that $X,Y,Z\in
{\cal L}(B)$ and $U,V,W\in {\cal L}(F)$. The connections $\nabla $, $%
{}^{B}\nabla $ and ${}^{F}\nabla $ are Levi-Civita connections on $M$, $B\,$%
and $F$, respectively. We will denote by $R,{}^{B}R$ and ${}^{F}R$ the
Riemann curvatures on $M$, $B\,$and $F$, respectively; $S$, ${}^{B}S$ and $%
{}^{F}S$ the Ricci tensors for the connections $\nabla $, ${}^{B}\nabla $
and ${}^{F}\nabla $, respectively; $r$,${}^{B}r$ and ${}^{F}r$ are the
scalar curvatures for the connections $\nabla $,${}^{B}\nabla $ and $%
{}^{F}\nabla $, respectively.

We need the following lemmas for later use. For more details see \cite%
{AKKY,Hat,Ponge,Unal-01}{\bf . }

\begin{lem}
\label{Lem1} Let $M=B_{f}\times {}_{h}F$ be a doubly warped product
manifold. If $X,Y\in $ ${\frak X}(B)$ and $V,W\in $ ${\frak X}(F)$, then
\[
{\rm \tan }\nabla _{X}Y=-\,\frac{{\rm grad}f}{f}g(X,Y),
\]%
\[
{\rm nor}\nabla _{X}Y\;{\rm is\;the\;lift\;of\;}\nabla _{X}Y\;{\rm on\;}B,
\]%
\[
\nabla _{X}V=\nabla _{V}X=\frac{Xh}{h}V+\frac{Vf}{f}X,
\]%
\[
{\rm nor}\nabla _{V}W=-\,\frac{{\rm grad}h}{h}g(V,W),
\]%
\[
{\rm tan}\nabla _{V}W\;{\rm is\;the\;lift\;of\;}\nabla _{V}W\;{\rm on\;}F,
\]%
where $\tan $ and $nor$ stand for tangent part to $F$ and normal part to $B$%
, respectively.
\end{lem}

\begin{lem}
\label{Lem2} Let $M=B_{f}\times {}_{h}F$ be a doubly warped product
manifold. If $X,Y,Z\in {\frak X}(B)$ and $U,V,W\in {\frak X}(F)$, then
\begin{eqnarray*}
R(X,Y)Z &=&{}^{B}\!R(X,Y)Z+\frac{\left\Vert {\rm grad\,}f\right\Vert ^{2}}{%
f^{2}}(g(X,Z)Y-g(Y,Z)X), \\
R(X,V)Y &=&\frac{H_{B}^{h}(X,Y)}{h}V+\frac{1}{f}g(X,Y)^{F}\nabla _{V}({\rm %
grad\,}f),
\end{eqnarray*}%
\begin{eqnarray*}
R(X,Y)V &=&\frac{Vf}{f}\left( \frac{Yh}{h}X-\frac{Xh}{h}Y\right) , \\
R(X,V)W &=&-\frac{H_{F}^{f}(V,W)}{f}X-\frac{1}{h}g(V,W)^{B}\nabla _{X}({\rm %
grad\,}h), \\
R(V,W)X &=&\frac{Xh}{h}\left( \frac{Wf}{f}V-\frac{Vf}{f}W\right) , \\
R(V,W)U &=&{}^{F}\!R(V,W)U+\frac{\left\Vert {\rm grad}h\right\Vert ^{2}}{%
h^{2}}(g(V,U)W-g(W,U)V),
\end{eqnarray*}%
where $H_{B}^{h}$ and $H_{F}^{f}$ are the Hessian of $h$ and $f$,
respectively.
\end{lem}

\begin{lem}
\label{Lem3} Let $M=B_{f}\times {}_{h}F$ be a doubly warped product
manifold. If $X,Y\in {\frak X}(B)$ and $V,W\in {\frak X}(F)$, then%
\begin{eqnarray*}
S(X,Y) &=&{}^{B}S(X,Y)-\frac{n_{2}}{h}H_{B}^{h}(X,Y) \\
&&-g(X,Y)\left( (n_{1}-1)\frac{\left\Vert {\rm grad}\,f\right\Vert ^{2}}{%
f^{2}}+\frac{\Delta _{F}f}{f}\right) , \\
S(X,V) &=&\frac{(n-2)(Xh)(Vf)}{hf}, \\
S(V,W) &=&{}^{F}S(V,W)-\frac{n_{1}}{f}H_{F}^{f}(V,W) \\
&&-g(V,W)\left( (n_{2}-1)\frac{\left\Vert {\rm grad\,}h\right\Vert ^{2}}{%
h^{2}}+\frac{\Delta _{B}h}{h}\right) ,
\end{eqnarray*}%
where $\Delta _{F}f$ and $\Delta _{B}h$ are the Laplacian of $f$ on $F$\ and
$h$ on $B$, respectively.
\end{lem}

\begin{lem}
\label{Lem4} Let $M=B_{f}\times {}_{h}F$ be a doubly warped product
manifold. Then
\begin{eqnarray}
r\ &=&\frac{1}{f^{2}}\ ^{B}r+\frac{1}{h^{2}}{}^{F}r-2n_{1}\frac{\Delta _{F}f%
}{f}-2n_{2}\frac{\Delta _{B}h}{h}  \nonumber \\
&&-n_{1}(n_{1}-1)\frac{\left\Vert {\rm grad\,}f\right\Vert ^{2}}{f^{2}}%
-n_{2}(n_{2}-1)\frac{\left\Vert {\rm grad\,}h\right\Vert ^{2}}{h^{2}}.
\label{scalarcurv}
\end{eqnarray}
\end{lem}

\section{Doubly Warped Product Manifolds Endowed with a Semi-symmetric
Metric Connection \label{Warped-ssmc}}

In this section, we consider doubly warped product manifolds with respect to
the semi-symmetric metric connection and find new expressions concerning
with curvature tensor, Ricci tensor and the scalar curvature admitting the
semi-symmetric metric connection where the associated vector field $P\in
{\frak X}(B)$ or $P\in {\frak X}(F)$. \medskip

We have the following results.

\begin{prop}
\label{Lemma1.1} Let $M=B_{f}\times {}_{h}F$ be a doubly warped product
manifold and let $\tilde{\nabla},{}^{B}\tilde{\nabla}$ and ${}^{F}\tilde{%
\nabla}$ be the semi-symmetric metric connections on $M,B$ and $F$,
respectively. If $X,Y\in {\frak X}(B),V,W\in {\frak X}(F)$ and $P\in {\frak X%
}(B)$, then%
\begin{equation}
{\rm \tan }\tilde{\nabla}_{X}Y=-\frac{1}{f}g(X,Y){\rm grad\,}f,  \label{eq-1}
\end{equation}%
\begin{equation}
{\rm nor}\tilde{\nabla}_{X}Y\;{\rm is\;the\;lift\;of\;\!}^{B}\tilde{\nabla}%
_{X}Y,  \label{eq-2}
\end{equation}%
\begin{equation}
\tilde{\nabla}_{X}V=\frac{Xh}{h}V+\frac{Vf}{f}X,  \label{eq-3}
\end{equation}%
\begin{equation}
\tilde{\nabla}_{V}X=\frac{Xh}{h}V+\frac{Vf}{f}X+\pi (X)V,  \label{eq-4}
\end{equation}%
\begin{equation}
{\rm nor}\tilde{\nabla}_{V}W=-\left( \frac{{\rm grad\,}h}{h}+P\right) g(V,W),
\label{eq-5}
\end{equation}%
\begin{equation}
{\rm tan}\tilde{\nabla}_{V}W\;{\rm is\;\;the\;\;lift\;\;of\;}^{F}\tilde{%
\nabla}_{V}W\;{\rm on\;}F,  \label{eq-6}
\end{equation}%
where $\tan $ and $nor$ stand for tangent part to $F$ and normal part to $B$%
, respectively.
\end{prop}

\noindent {\bf Proof. }Consider{\bf \ }$P\in {\frak X}(B)$ and using (\ref%
{semsymmetric}), first equation of Lemma \ref{Lem1}, we get (\ref{eq-1}).
Simailarly, we can easily find the others results by using (\ref%
{semsymmetric}) and Lemma \ref{Lem1}. $\blacksquare $

\begin{prop}
\label{Lemma1.2} Let $M=B_{f}\times {}_{h}F$ be a doubly warped product
manifold and let $\tilde{\nabla},{}^{B}\tilde{\nabla}$ and ${}^{F}\tilde{%
\nabla}$ be the semi-symmetric metric connections on $M,B$ and $F,$
respectively. If $X,Y\in {\frak X}(B),V,W\in {\frak X}(F)$ and $P\in {\frak X%
}(F)$, then%
\begin{equation}
{\rm \tan }\tilde{\nabla}_{X}Y=-\left( \frac{{\rm grad\,}f}{f}+P\right)
g(X,Y),  \label{eq-7}
\end{equation}%
\begin{equation}
{\rm nor}\tilde{\nabla}_{X}Y\;{\rm is\;the\;lift\;of}\!\;^{B}\tilde{\nabla}%
_{X}Y,  \label{eq-8}
\end{equation}%
\begin{equation}
\tilde{\nabla}_{X}V=\frac{Xh}{h}V+\frac{Vf}{f}X+\pi (V)X,  \label{eq-9}
\end{equation}%
\begin{equation}
\tilde{\nabla}_{V}X=\frac{Xh}{h}V+\frac{Vf}{f}X,  \label{eq-10}
\end{equation}%
\begin{equation}
{\rm nor}\tilde{\nabla}_{V}W=-\frac{{\rm grad\,}h}{h}g(V,W),  \label{eq-11}
\end{equation}%
\begin{equation}
{\rm tan}\tilde{\nabla}_{V}W\;{\rm is\;the\;lift\;of}\!\;^{F}\tilde{\nabla}%
_{V}W\;{\rm on\;}F.  \label{eq-12}
\end{equation}
\end{prop}

\noindent {\bf Proof.} Consider{\bf \ }$P\in {\frak X}(F)$ and using (\ref%
{semsymmetric}), first equation of Lemma \ref{Lem1}, we get (\ref{eq-7}).
Simailarly, we can easily find the others results by using (\ref%
{semsymmetric}) and Lemma \ref{Lem1}.. $\blacksquare $

\begin{prop}
\label{Lemma2.1} Let $M=B_{f}\times {}_{h}F$ be a doubly warped product
manifold and let $\tilde{R},{}^{B}\tilde{R}$ and ${}^{F}\tilde{R}$ be the
Riemannian curvature tensors with respect to the semi-symmetric metric
connections $\tilde{\nabla},{}^{B}\tilde{\nabla}$ and ${}^{F}\tilde{\nabla}$%
, respectively. If $X,Y,Z\in {\frak X}(B),U,V,W\in {\frak X}(F)$ and $P\in
{\frak X}(B)$, then%
\begin{eqnarray*}
\tilde{R}(X,Y)Z &=&\!^{B}\tilde{R}(X,Y)Z+\frac{\left\Vert {\rm grad\,}%
f\right\Vert ^{2}}{f^{2}}(g(X,Z)Y-g(Y,Z)X) \\
&&+(g(Y,Z)\pi (X)-g(X,Z)\pi (Y))\frac{{\rm grad\,}f}{f}, \\
\tilde{R}(V,X)Y &=&\left( \frac{H_{B}^{h}(X,Y)}{h}{\bf +}\pi (X)\pi
(Y)-g(Y,^{B}\nabla _{X}P)\right) V+\frac{Vf}{f}\pi (Y)X \\
&&-\left( \frac{Vf}{f}P+\frac{Ph}{h}V+\pi (P)V-\frac{1}{f}^{F}\nabla _{V}%
{\rm grad\,}f\right) g(X,Y), \\
\tilde{R}(X,Y)V &=&\left( \frac{(Vf)(Yh)}{hf}+\frac{(Vf)}{f}\pi (Y)\right) X
\\
&&-\left( \frac{(Vf)(Xh)}{hf}+\frac{(Vf)}{f}\pi (X)\right) Y, \\
\tilde{R}(V,W)X &=&\left( \frac{(Wf)(Xh)}{hf}-\frac{(Wf)}{f}\pi (X)\right) V
\\
&&-\left( \frac{(Vf)(Xh)}{hf}-\frac{(Vf)}{f}\pi (X)\right) W, \\
\tilde{R}(X,V)W &=&-\frac{H_{F}^{f}(V,W)}{f}X-\frac{Wf}{f}\pi (X)V \\
&&-g(V,W)\left( \frac{^{B}\nabla _{X}{\rm grad}h}{h}+\frac{Ph}{h}%
X+\!^{B}\nabla _{X}P\right. \\
&&\left. +\pi (P)X-\pi (X)P-\frac{{\rm grad}f}{f}\pi (X)\right) , \\
\tilde{R}(U,V)W &=&\!^{F}R(U,V)W-\frac{Uf}{f}g(V,W)P+\frac{Vf}{f}g(U,W)P \\
&&-\left( \frac{\left\Vert {\rm grad}h\right\Vert ^{2}}{h^{2}}+\frac{2Ph}{h}%
+\pi (P)\right) (g(V,W)U-g(U,W)V.
\end{eqnarray*}
\end{prop}

\noindent {\bf Proof.} Assume that $M=B_{f}\times {}_{h}F$ is a doubly
warped product and $R$ and $\tilde{R}$ denote the curvature tensors with
respect to the Levi-Civita connection and the semi-symmetric metric
connection, respectively.

In view of the (\ref{curvature}), Lemma \ref{Lem1} and Lemma \ref{Lem2}, we
can write
\begin{eqnarray*}
\tilde{R}(X,Y)Z &=&\!^{B}R(X,Y)Z+\frac{\left\Vert {\rm grad\,}f\right\Vert
^{2}}{f^{2}}(g(X,Z)Y-g(Y,Z)X) \\
&&+g(Z,^{B}\nabla _{X}P)Y-g(Z,^{B}\nabla _{Y}P)X \\
&&+g(X,Z)(^{B}\nabla _{Y}P+\pi (P)Y-\pi (Y)P) \\
&&-g(Y,Z)(^{B}\nabla _{X}P+\pi (P)X-\pi (X)P) \\
&&+(g(Y,Z)\pi (X)-g(X,Z)\pi (Y))\frac{{\rm grad\,}f}{f} \\
&&+\pi (Z)(\pi (Y)X-\pi (X)Y).
\end{eqnarray*}%
By using equation (\ref{curvature}), we can write%
\begin{eqnarray}
\tilde{R}(V,X)Y &=&R(V,X)Y+g(Y,\nabla _{V}P)X-g(Y,\nabla _{X}P)V  \nonumber
\\
&&-g(X,Y)[\nabla _{V}P+\pi (P)V-\pi (V)P]  \label{RVXY} \\
&&+\pi (Y)[\pi (X)V-\pi (V)X].  \nonumber
\end{eqnarray}%
Since $P\in {\frak X}(B)$ and by making use of Lemma \ref{Lem1} and Lemma %
\ref{Lem2}, we get%
\begin{eqnarray*}
\tilde{R}(V,X)Y &=&\left( \frac{H_{B}^{h}(X,Y)}{h}{\bf +}\pi (X)\pi
(Y)-g(Y,^{B}\nabla _{X}P)\right) V+\frac{Vf}{f}\pi (Y)X \\
&&-\left( \frac{Vf}{f}P+\frac{Ph}{h}V+\pi (P)V-\frac{1}{f}^{F}\nabla
_{V}grad\,f\right) g(X,Y).
\end{eqnarray*}%
Replacing $Z$ by $V$ in equation (\ref{curvature}), we get%
\begin{equation}
\tilde{R}(X,Y)V=R(X,Y)V+g(V,\nabla _{X}P)Y-g(V,\nabla _{Y}P)X.  \label{RXYV}
\end{equation}%
Using Lemma \ref{Lem1} and Lemma \ref{Lem2}, we get
\begin{eqnarray*}
\tilde{R}(X,Y)V &=&\left( \frac{(Vf)(Yh)}{hf}+\frac{(Vf)}{f}\pi (Y)\right) X
\\
&&-\left( \frac{(Vf)(Xh)}{hf}+\frac{(Vf)}{f}\pi (X)\right) Y.
\end{eqnarray*}%
By making use of (\ref{curvature}) and Lemma \ref{Lem1} Lemma \ref{Lem2}, we
get
\begin{eqnarray}
\tilde{R}(V,W)X &=&\left( \frac{(Wf)(Xh)}{hf}-\frac{(Wf)}{f}\pi (X)\right) V
\nonumber \\
&&-\left( \frac{(Vf)(Xh)}{hf}-\frac{(Vf)}{f}\pi (X)\right) W.
\label{eq-RVWX}
\end{eqnarray}%
From the equation (\ref{curvature}), we find
\begin{eqnarray}
\tilde{R}(X,V)W &=&R(X,V)W+g(W,\nabla _{X}P)V-g(W,\nabla _{V}P)X  \nonumber
\\
&&-g(V,W)(\nabla _{X}P+\pi (P)X-\pi (X)P).  \label{RXVW}
\end{eqnarray}
Using Lemma \ref{Lem1} and Lemma \ref{Lem2} in (\ref{RXVW}), we have%
\begin{eqnarray*}
\tilde{R}(X,V)W &=&-\frac{H_{F}^{f}(V,W)}{f}X-\frac{Wf}{f}\pi (X)V \\
&&-g(V,W)\left( \frac{^{B}\nabla _{X}{\rm grad}h}{h}+\frac{Ph}{h}%
X+\!^{B}\nabla _{X}P\right. \\
&&\left. +\pi (P)X-\pi (X)P-\frac{{\rm grad}f}{f}\pi (X)\right) .
\end{eqnarray*}%
In view of the equation (\ref{curvature}), we have%
\begin{eqnarray}
\tilde{R}(U,V)W &=&R(U,V)W+g(W,\nabla _{U}P)V-g(W,\nabla _{V}P)U  \nonumber
\\
&&+g(U,W)\nabla _{V}P-g(V,W)\nabla _{U}P  \nonumber \\
&&+\pi (P)[g(U,W)V-g(V,W)U].  \label{RUVW}
\end{eqnarray}%
By making use of Lemma \ref{Lem1} and Lemma \ref{Lem2} in (\ref{RUVW}), we
obtain
\begin{eqnarray*}
\tilde{R}(U,V)W &=&\!^{F}R(U,V)W-\frac{Uf}{f}g(V,W)P+\frac{Vf}{f}g(U,W)P \\
&&-\left( \frac{\left\Vert {\rm grad}h\right\Vert ^{2}}{h^{2}}+\frac{2Ph}{h}%
+\pi (P)\right) (g(V,W)U-g(U,W)V.
\end{eqnarray*}%
Hence, the proof is completed. $\blacksquare $

\begin{prop}
\label{Lemma2.2} Let $M=B_{f}\times {}_{h}F$ be a doubly warped product
manifold and let $\tilde{R},{}^{B}\tilde{R}$ and ${}^{F}\tilde{R}$ be the
Riemannian curvature tensors with respect to the semi-symmetric metric
connections $\tilde{\nabla},{}^{B}\tilde{\nabla}$ and ${}^{F}\tilde{\nabla}$%
, respectively. If $X,Y,Z\in {\frak X}(B),U,V,W\in {\frak X}(F)$ and $P\in
{\frak X}(F)$, then%
\begin{eqnarray*}
\tilde{R}(X,Y)Z &=&\!^{B}R(X,Y)Z+\left( g(X,Z)\frac{Yh}{h}-g(Y,Z)\frac{Xh}{h}%
\right) \left( P-\frac{{\rm gradf}}{f}\right) \\
&&+\left( \frac{\left\Vert {\rm grad}f\right\Vert ^{2}}{f^{2}}+\frac{2Pf}{f}%
+\pi (P)\right) (g(X,Z)Y-g(Y,Z)X), \\
\tilde{R}(V,X)Y &=&-\left( \frac{H_{B}^{h}(X,Y)}{h}+\frac{Pf}{f}g(X,Y)+\pi
(P)g(X,Y)\right) V-\frac{(Yh)}{h}\pi (V)X \\
&&+g(X,Y)\left( \frac{1}{h}\pi (V){\rm grad}h+\pi (V)P-\,^{F}\nabla _{V}P-%
\frac{1}{f}^{F}\nabla _{V}{\rm grad}f\right) , \\
\tilde{R}(X,Y)V &=&\left( \frac{(Vf)(Yh)}{hf}-\frac{(Yh)}{h}\pi (V)\right) X
\\
&&-\left( \frac{(Vf)(Xh)}{hf}-\frac{(Xh)}{h}\pi (V)\right) Y, \\
\tilde{R}(V,W)X &=&\left( \frac{(Wf)(Xh)}{hf}+\frac{(Xh)}{h}\pi (W)\right) V
\\
&&-\left( \frac{(Vf)(Xh)}{hf}+\frac{(Xh)}{h}\pi (V)\right) W, \\
\tilde{R}(X,V)W &=&-g(V,W)\left( \frac{1}{h}^{B}\nabla _{X}{\rm grad}h+\frac{%
Pf}{f}X+\pi (P)X+\frac{(Xh)}{h}P\right) \\
&&-\left( \frac{H_{F}^{f}(V,W)}{f}-\pi (V)\pi (W)+g(W,^{F}\nabla
_{V}P)\right) X+\frac{(Xh)}{h}\pi (W)V, \\
\tilde{R}(U,V)W &=&\!^{F}\tilde{R}(U,V)W \\
&&-\frac{\left\Vert {\rm grad}h\right\Vert ^{2}}{h^{2}}(g(V,W)U-g(U,W)V) \\
&&+\left( \pi (U)g(V,W)-g(U,W)\pi (V)\right) \frac{{\rm grad}h}{h}.
\end{eqnarray*}
\end{prop}

\noindent {\bf Proof.} Assume that the associated vector field $P\in {\frak X%
}(F)$. Then the equation (\ref{curvature}) can be written as
\begin{eqnarray*}
\tilde{R}(X,Y)Z &=&R(X,Y)Z+g(Z,\nabla _{X}P)Y-g(Z,\nabla _{Y}P)X \\
&&+g(X,Z)\nabla _{Y}P-g(Y,Z)\nabla _{X}P \\
&&+\pi (P)(g(X,Z)Y-g(Y,Z)X).
\end{eqnarray*}%
By the use of Lemma \ref{Lem1} and Lemma \ref{Lem2}, the above equation
gives us%
\begin{eqnarray*}
\tilde{R}(X,Y)Z &=&\!^{B\,}R(X,Y)Z \\
&&+\left( \frac{\left\Vert {\rm grad\,}f\right\Vert ^{2}}{f^{2}}+\frac{2Pf}{f%
}+\pi (P)\right) (g(X,Z)Y-g(Y,Z)X) \\
&&+(g(Y,Z)(Xh)-g(X,Z)(Yh))\left( \frac{{\rm grad\,}f}{hf}-\frac{P}{h}\right)
.
\end{eqnarray*}%
By (\ref{RVXY}), Lemma \ref{Lem1} and Lemma \ref{Lem2}, we obtain
\begin{eqnarray*}
\tilde{R}(V,X)Y &=&-\left( \frac{H_{B}^{h}(X,Y)}{h}+\frac{Pf}{f}g(X,Y)+\pi
(P)g(X,Y)\right) V-\frac{(Yh)}{h}\pi (V)X \\
&&+g(X,Y)\left( \frac{1}{h}\pi (V){\rm grad}h+\pi (V)P-\,^{F}\nabla _{V}P-%
\frac{1}{f}^{F}\nabla _{V}{\rm grad}f\right) ,
\end{eqnarray*}%
Replacing $Z$ by $V$ in equation (\ref{curvature}), we get
\[
\tilde{R}(X,Y)V=R(X,Y)V+g(V,\nabla _{X}P)Y-g(V,\nabla _{Y}P)X.
\]%
Using Lemma \ref{Lem1} and Lemma \ref{Lem2}, we obtain
\begin{eqnarray*}
\tilde{R}(X,Y)V &=&\left( \frac{(Vf)(Yh)}{hf}-\frac{Yh}{h}\pi (V)\right) X \\
&&-\left( \frac{(Vf)(Xh)}{hf}-\frac{Xh}{h}\pi (V)\right) Y.
\end{eqnarray*}%
From equation (\ref{curvature}), we get
\[
\tilde{R}(V,W)X=R(V,W)X+g(X,\nabla _{V}P)W-g(X,\nabla _{W}P)V.
\]%
By Lemma \ref{Lem1} and Lemma \ref{Lem2}, we get
\begin{eqnarray*}
\tilde{R}(V,W)X &=&\left( \frac{(Wf)(Xh)}{hf}+\frac{(Xh)}{h}\pi (W)\right) V
\\
&&-\left( \frac{(Vf)(Xh)}{hf}-\frac{(Xh)}{h}\pi (V)\right) W.
\end{eqnarray*}%
From equation (\ref{curvature}), we get
\begin{eqnarray*}
\tilde{R}(X,V)W &=&R(X,V)W+g(W,\nabla _{X}P)V-g(W,\nabla _{V}P)X \\
&&-g(V,W)\nabla _{X}P-\pi (P)g(V,W)X+\pi (V)\pi (W)X.
\end{eqnarray*}%
Using Lemma \ref{Lem1} and Lemma \ref{Lem2}, we get
\begin{eqnarray*}
\tilde{R}(X,V)W &=&-g(V,W)\left( \frac{1}{h}^{B}\nabla _{X}{\rm grad}h+\frac{%
Pf}{f}X+\pi (P)X+\frac{(Xh)}{h}P\right) \\
&&-\left( \frac{H_{F}^{f}(V,W)}{f}-\pi (V)\pi (W)+g(W,^{F}\nabla
_{V}P)\right) X+\frac{(Xh)}{h}\pi (W)V.
\end{eqnarray*}%
From equation (\ref{curvature}), we have
\begin{eqnarray}
\tilde{R}(U,V)W &=&R(U,V)W+g(W,\nabla _{U}P)V-g(W,\nabla _{V}P)U  \nonumber
\\
&&+g(U,W)\nabla _{V}P-g(V,W)\nabla _{U}P  \nonumber \\
&&+\pi (P)(g(U,W)V-g(V,W)U)  \nonumber \\
&&+(g(V,W)\pi (U)-g(U,W)\pi (U))P  \nonumber \\
&&+\pi (W)(\pi (V)U-\pi (U)V).  \label{RUVW1}
\end{eqnarray}%
By use of Lemma \ref{Lem1} and Lemma \ref{Lem2}in above equation, we obtain
\begin{eqnarray*}
\tilde{R}(U,V)W &=&\!^{F}\tilde{R}(U,V)W \\
&&-\frac{\left\Vert {\rm grad}h\right\Vert ^{2}}{h^{2}}(g(V,W)U-g(U,W)V) \\
&&+\left( \pi (U)g(V,W)-g(U,W)\pi (V)\right) \frac{{\rm grad}h}{h}.
\end{eqnarray*}%
Thus, we complete the proof. $\blacksquare $

As a consequence of Proposition \ref{Lemma2.1} and Proposition \ref{Lemma2.2}%
, by a contraction of the curvature tensors we obtain the Ricci tensors of
the doubly warped product with respect to the semi-symmetric metric
connection as follows:

\begin{cor}
\label{Corollary1.1} Let $M=B_{f}\times {}_{h}F$ be a doubly warped product
manifold and let $\tilde{S},{}^{B}\tilde{S}$ and ${}^{F}\tilde{S}$ be the
Ricci tensors with respect to the semi-symmetric metric connections $\tilde{%
\nabla},{}^{B}\tilde{\nabla}$ and ${}^{F}\tilde{\nabla}$, respectively. If $%
X,Y\in {\frak X}(B),V,W\in {\frak X}(F)$ and $P\in {\frak X}(B)$, then%
\begin{eqnarray*}
\tilde{S}(X,Y) &=&\!^{B}\tilde{S}(X,Y)-\frac{n_{2}}{h}H_{B}^{h}(X,Y)+n_{2}%
\pi (X)\pi (Y)-n_{2}g(Y,^{B}\nabla _{X}P) \\
&&-\left( (n_{1}-1)\frac{\left\Vert gradf\right\Vert ^{2}}{f^{2}}+n_{2}\pi
(P)+n_{2}\frac{Ph}{h}+\frac{1}{f}\triangle f\right) g(X,Y), \\
\tilde{S}(X,V) &=&(n_{1}-1)\frac{(Vf)(Xh)}{hf}+(n-2)\frac{(Vf)}{f}\pi (X), \\
\tilde{S}(V,X) &=&(n_{2}-1)\frac{(Vf)(Xh)}{hf}-(n-2)\frac{(Vf)}{f}\pi (X), \\
\tilde{S}(V,W) &=&^{F}S(V,W)-\frac{n_{1}}{f}H_{F}^{f}(V,W) \\
&&-\left( {\rm div}P+\frac{1}{h}\triangle h+n_{1}\frac{Ph}{h}+(n_{1}-1)\pi
(P)\right) g(V,W) \\
&&-(n_{2}-1)\left( \frac{\left\Vert {\rm grad}h\right\Vert ^{2}}{h^{2}}+%
\frac{2Ph}{h}+\pi (P)\right) g(V,W).
\end{eqnarray*}
\end{cor}

\begin{cor}
\label{Corollary1.2} Let $M=B_{f}\times {}_{h}F$ be a doubly warped product
manifold and let $\tilde{S},{}^{B}\tilde{S}$ and ${}^{F}\tilde{S}$ be the
Ricci tensors with respect to the semi-symmetric metric connections $\tilde{%
\nabla},{}^{B}\tilde{\nabla}$ and ${}^{F}\tilde{\nabla}$, respectively. If $%
X,Y\in {\frak X}(B),V,W\in {\frak X}(F)$ and $P\in {\frak X}(F)$, then%
\begin{eqnarray*}
\tilde{S}(X,Y) &=&\!^{B}S(X,Y)-n_{2}\frac{H_{B}^{h}(X,Y)}{h} \\
&&-(n_{1}-1)\left( \frac{\left\Vert gradf\right\Vert ^{2}}{f^{2}}+\frac{2Pf}{%
f}\right) g(X,Y) \\
&&-\left( \frac{1}{f}\!^{F}\bigtriangleup f+(n-2)\pi (P)+n_{2}\frac{Pf}{f}+%
{\rm div}P\right) g(X,Y), \\
\tilde{S}(X,V) &=&(n_{1}-1)\frac{(Vf)(Xh)}{hf}-(n-2)\frac{Xh}{h}\pi (V), \\
\tilde{S}(V,W) &=&\!^{F}\tilde{S}(V,W) \\
&&-n_{1}g(W,^{F}\nabla _{V}P)-\frac{n_{1}}{f}H_{F}^{f}(V,W)+n_{1}\pi (V)\pi
(W) \\
&&-\left( \frac{1}{h}^{B}\bigtriangleup h+n_{1}\frac{Pf}{f}+(n_{2}-1)\frac{%
\left\Vert {\rm grad}h\right\Vert ^{2}}{h^{2}}+n_{1}\pi (P)\right) g(V,W).
\end{eqnarray*}
\end{cor}

As a consequence of Corollary \ref{Corollary1.1} and Corollary \ref%
{Corollary1.2}, by a contraction of the Ricci tensors we get scalar
curvatures of the doubly warped product with respect to the semi-symmetric
metric connection as follows:

\begin{cor}
\label{Corollary2.1} Let $M=B_{f}\times {}_{h}F$ be a doubly warped product
manifold and $P\in {\frak X}(B)$. Let $\tilde{r},{}^{B}\tilde{r}$ and${}^{F}%
\tilde{r}$ be the scalar curvatures with respect to the semi-symmetric
metric connections $\tilde{\nabla},{}^{B}\tilde{\nabla}$ and ${}^{F}\tilde{%
\nabla}$, respectively. Then
\begin{eqnarray*}
\tilde{r} &=&\!\frac{^{B}\tilde{r}}{f^{2}}+\!\frac{^{F}r}{h^{2}}%
-n_{1}(n_{1}-1)\frac{\left\Vert gradf\right\Vert ^{2}}{f^{2}}-n_{2}(n_{2}-1)%
\frac{\left\Vert {\rm grad}h\right\Vert ^{2}}{h^{2}} \\
&&-2n_{2}(n-1)\frac{Ph}{h}-2n_{1}\frac{^{F}\Delta f}{f}-2n_{2}\frac{%
^{B}\Delta h}{h}-2n_{2}{\rm div}P \\
&&-n_{2}(n+n_{1}-3)\pi (P).
\end{eqnarray*}
\end{cor}

\begin{cor}
\label{Corollary2.2} Let $M=B_{f}\times {}_{h}F$ be a doubly warped product
and $P\in {\frak X}(F)$. Let $\tilde{r},{}^{B}\tilde{r}$ and ${}^{F}\tilde{r}
$ be the scalar curvatures with respect to the semi-symmetric metric
connections $\tilde{\nabla},{}^{B}\tilde{\nabla}$ and ${}^{F}\tilde{\nabla}$%
, respectively. Then
\begin{eqnarray*}
\tilde{r} &=&\!\frac{^{B}r}{f^{2}}+\!\frac{^{F}\tilde{r}}{h^{2}}%
-n_{1}(n_{1}-1)\frac{\left\Vert {\rm grad}f\right\Vert ^{2}}{f^{2}}%
-n_{2}(n_{2}-1)\frac{\left\Vert {\rm grad}h\right\Vert ^{2}}{h^{2}} \\
&&-2n_{1}(n-1)\frac{Pf}{f}-2n_{1}\frac{^{F}\Delta f}{f}-2n_{2}\frac{%
^{B}\Delta h}{h}-2n_{1}{\rm div}P \\
&&-n_{1}(n+n_{2}-3)\pi (P).
\end{eqnarray*}
\end{cor}

\begin{rem-new}
Doubly warped product manifolds with the semi-symmetric metric connection
has been also studied by Sular \cite{Sular}, where the author has taken the
associated vector field $P\in {\frak X}(M)$ as $P=P_{B}+P_{F}$, $P_{B}$ and $%
P_{F}$ are the components of $P$ on $B$ and $F$, respectively.
\end{rem-new}

\section{Einstein doubly warped product manifolds endowed with the
semi-symmetric metric connection}

In this section, we consider Einstein doubly warped products and
Einstein-like doubly warped product of class ${\cal A}$ endowed with the
semi-symmetric metric connection.

\begin{th}
\label{Corollary1.1 copy(1)} Let $M=B_{f}\times {}_{h}F$ be a doubly warped
product manifold $(n>2)$ and $P\in {\frak X}(B)$. Then $(M,\tilde{\nabla})$
is an Einstein manifold with Einstein constant $\mu $ if and only if
\begin{eqnarray}
^{B}\tilde{S}(X,Y) &=&\!\frac{n_{2}}{h}H_{B}^{h}(X,Y)-n_{2}\pi (X)\pi
(Y)+n_{2}g(Y,^{B}\nabla _{X}P)  \nonumber \\
&&+\left( (n_{1}-1)\frac{\left\Vert gradf\right\Vert ^{2}}{f^{2}}+n_{2}\pi
(P)+n_{2}\frac{Ph}{h}+\frac{1}{f}\triangle f+\mu \right) g(X,Y),
\label{eq-e1}
\end{eqnarray}%
\begin{eqnarray}
^{F}S(V,W) &=&\frac{n_{1}}{f}H_{F}^{f}(V,W)  \nonumber \\
&&+\left( {\rm div}P+\frac{1}{h}\triangle h+n_{1}\frac{Ph}{h}+(n_{1}-1)\pi
(P)\right) g(V,W)  \nonumber \\
&&+(n_{2}-1)\left( \frac{\left\Vert {\rm grad}h\right\Vert ^{2}}{h^{2}}+%
\frac{2Ph}{h}+\pi (P)+\mu \right) g(V,W).  \label{eq-e3}
\end{eqnarray}
\end{th}

\noindent {\bf Proof.} The proof follows from Corollary \ref{Corollary1.1}
and Corollary \ref{Corollary1.2}. $\blacksquare $

\begin{th}
Let $(M,\tilde{\nabla})$ be an Einstein doubly warped product manifold with
Einstein constant $\mu $ and $P\in {\frak X}(B)$. Suppose that $B$ is
compact Riemannian manifold, $F$ is complete Riemannian manifold, $n_{1}$ , $%
n_{2}$ $\geq 2$, and $H_{F}^{f}(V,W)$ is a constant multiple of $g_{F}$. If $%
\mu \leq 0$, $Ph\geq 0$, $\pi (P)\leq 0$, ${\rm div}P\leq 0$ then $M$ is a
warped product manifold.
\end{th}

\noindent {\bf Proof.} Since $H_{F}^{f}(V,W)$ is a constant multiple of $%
g_{F}$, so by using the result of Tashiro \cite[Theorem 2]{Tashiro}, we can
say that $F$ is a Euclidean space, then Ricci tensor of $F$ is zero. By (\ref%
{eq-e3}), we get
\begin{eqnarray*}
0 &=&\left( \frac{cn_{1}}{f}+(n_{2}-1)\left( \left\Vert {\rm grad}%
h\right\Vert ^{2}+2hPh+h^{2}\pi (P)+\mu h^{2}\right) \right.  \\
&&+\left. h^{2}{\rm div}P+h\triangle h+n_{1}hPh+(n_{1}-1)h^{2}\pi (P)\right)
g_{F}(V,W).
\end{eqnarray*}%
So
\begin{eqnarray}
0 &=&\left( \frac{cn_{1}}{f}+(n_{2}-1)\left( \left\Vert {\rm grad}%
h\right\Vert ^{2}+2hPh+h^{2}\pi (P)+\mu h^{2}\right) \right.   \nonumber \\
&&+\left. h^{2}{\rm div}P+h\triangle h+n_{1}hPh+(n_{1}-1)h^{2}\pi (P)\right)
,  \label{eq-min}
\end{eqnarray}%
Let $x\in B$ such that $h(x)$ is maximum of $h$ on $B$. Therefore ${\rm grad}%
h(x)=0$ and $\triangle h(x)\leq 0$. Then $Ph(x)=g({\rm grad}h(x),P)=0$. So
the eq(\ref{eq-min}) at the point $x$ is
\begin{eqnarray}
0 &=&\left( (n_{2}-1)\left( h^{2}(x)\pi (P)+\mu h^{2}(x)\right) \right.
\nonumber \\
&&+\left. h^{2}(x){\rm div}P+h(x)\triangle h(x)+(n_{1}-1)h^{2}(x)\pi
(P)\right) .  \label{eq-min-1}
\end{eqnarray}%
By eq(\ref{eq-min}) and eq(\ref{eq-min-1}), we have
\begin{eqnarray}
0 &=&\left( \frac{cn_{1}}{f}+(n_{2}-1)\left( \left\Vert {\rm grad}%
h\right\Vert ^{2}+2hPh+(h^{2}-h^{2}(x))\pi (P)+\mu (h^{2}-h^{2}(x))\right)
\right.   \nonumber \\
&&+\left. (h^{2}-h^{2}(x)){\rm div}P+h\triangle h-h(x)\triangle
h(x)+n_{1}hPh\right.   \nonumber \\
&&+\left. (n_{1}-1)(h^{2}-h^{2}(x))\pi (P)\right) .  \label{eq-min-2}
\end{eqnarray}%
Since $\mu \leq 0$, $Ph\geq 0$, $\pi (P)\leq 0$, ${\rm div}P\leq 0$, so
equation (\ref{eq-min-2}) implies that $h\triangle h\leq 0$, which shows
that Laplacian has constant sign and hence $h$ is constant. $\blacksquare $

\begin{th}
Let $(M,g)$ be an Einstein doubly warped product manifold $M=I_{f}\times
{}_{h}F$ with respect to the semi-symmetric metric connection, where $I$ is
an open interval, {\rm dim}$I=1$ and {\rm dim}$F=n-1(n\geq 3)$. If $P\in
{\frak X}(I)$, then $f$ is constant on $F$ or $Ph=h\pi (P)$.
\end{th}

\noindent {\bf Proof.} Let $(M,g)$ be a doubly warped product $M=I_{f}\times
{}_{h}F$, where $I$ is an open interval, {\rm dim}$I=1$ and {\rm dim}$%
F=n-1(n\geq 3)$. By Corollary \ref{Corollary1.1}, we have
\[
\tilde{S}(X,V)=(n-2)\frac{Vf}{f}\pi (X),
\]%
\begin{equation}
\tilde{S}(V,X)=(n-2)\frac{(Vf)}{f}\left( \frac{Xh}{h}-\pi (X)\right) .
\label{eq-dwp-semi-4}
\end{equation}%
Since $M$ is an Einstein manifold with respect to the semi-symmetric metric
connection, we can write
\[
\tilde{S}(P,V)=\alpha g(P,V),
\]%
where $\alpha $ is constant. But $g(P,V)=0$, therefore $\tilde{S}(P,V)=0$.
By (\ref{eq-dwp-semi-4}), we have either $Vf=0$ or $\displaystyle\frac{Ph}{h}%
=\pi (P)$. Therefore either $f$ is constant or $Ph=h\pi (P)$. $\blacksquare $

\begin{defn-new}
\cite{Sumitomo} A Riemannian manifold $(M,g)$ is said to admit a
cyclic-Ricci parallel tensor or is Einstein-like of class ${\cal A}$ if%
\[
(\nabla _{X}S)(Y,Z)+(\nabla _{Y}S)(Z,X)+(\nabla _{Z}S)(X,Y)=0
\]%
for any vector fields $X,Y,Z\in {\frak X}(M)$ or equivalently $(\nabla
_{X}S)(X,X)=0$.
\end{defn-new}

\begin{prop}
Let $M=B_{f}\times {}_{h}F$ be an Einstein-like doubly warped product
manifold of class ${\cal A}$ with respect to the semi-symmetric metric
connection $\tilde{\nabla}$ and $P\in {\frak X}(B)$. The Riemannian manifold
$B$ is an Einstein-like manifold of class ${\cal A}$ with respect to the
semi-symmetric metric connection $\tilde{\nabla}$ if and only if $^{B}\nabla
_{X}Y=\frac{Xh}{2h}Y$, $^{B}\nabla _{X}P=0$ and
\begin{eqnarray*}
0 &=&n_{2}\pi (X)\left( \pi (X)\frac{Xh}{h}-\frac{\left\Vert {\rm grad}%
h\right\Vert ^{2}}{2h^{2}}g(X,X)\right) \\
&&+\frac{Xh}{h}g(X,X)\left( (n_{2}-1)\frac{\left\Vert gradf\right\Vert ^{2}}{%
f^{2}}-n_{2}\pi (P)-n_{2}\frac{1}{f}\triangle f\right) .
\end{eqnarray*}
\end{prop}

\noindent {\bf Proof. }By using (\ref{eq-1}), (\ref{eq-2}) and Corollary\ref%
{Corollary1.1}, we have
\begin{eqnarray*}
(\tilde{\nabla}_{X}\tilde{S})(X,X) &=&X\tilde{S}(X,X)-\tilde{S}(\tilde{\nabla%
}_{X}X,X)-\tilde{S}(X,\tilde{\nabla}_{X}X) \\
&=&X\left( ^{B}\tilde{S}(X,X)-\frac{n_{2}}{h}H_{B}^{h}(X,X)+n_{2}\pi (X)\pi
(X)-n_{2}g(X,^{B}\nabla _{X}P)\right) \\
&&-X\left( \left( (n_{1}-1)\frac{\left\Vert gradf\right\Vert ^{2}}{f^{2}}%
+n_{2}\pi (P)+n_{2}\frac{Ph}{h}+\frac{1}{f}\triangle f\right) g(X,X)\right)
\\
&&-\tilde{S}(\tilde{\nabla}_{X}^{B}X-\frac{1}{f}g(X,X){\rm grad}f,X)-\tilde{S%
}(X,\tilde{\nabla}_{X}^{B}X-\frac{1}{f}g(X,X){\rm grad}f),
\end{eqnarray*}%
\begin{eqnarray*}
(\tilde{\nabla}_{X}\tilde{S})(X,X) &=&(\tilde{\nabla}_{X}^{B}\tilde{S})(X,X)+%
\frac{n_{2}}{h^{2}}XhH_{B}^{h}(X,X)-\frac{n_{2}}{h}X\left(
H_{B}^{h}(X,X)\right) \\
&&+n_{2}X(\pi (X))^{2}-n_{2}Xg(X,^{B}\nabla _{X}P) \\
&&-\left( (n_{1}-1)\frac{\left\Vert gradf\right\Vert ^{2}}{f^{2}}+n_{2}\pi
(P)+n_{2}\frac{Ph}{h}+\frac{1}{f}\triangle f\right) Xg(X,X) \\
&&-\left( (n_{1}-1)X\frac{\left\Vert gradf\right\Vert ^{2}}{f^{2}}+n_{2}X\pi
(P)+n_{2}X\frac{Ph}{h}+X\frac{1}{f}\triangle f\right) g(X,X) \\
&&+\frac{1}{f}g(X,X)\tilde{S}({\rm grad}f,X)+\frac{1}{f}g(X,X)\tilde{S}(X,%
{\rm grad}f).
\end{eqnarray*}%
Since $^{B}\nabla _{X}Y=\frac{Xh}{2h}Y$ and $^{B}\nabla _{X}P=0$, we have
\begin{eqnarray*}
(\tilde{\nabla}_{X}\tilde{S})(X,X) &=&(\tilde{\nabla}_{X}^{B}\tilde{S}%
)(X,X)+n_{2}\pi (X)\left( \pi (X)\frac{Xh}{h}-\frac{\left\Vert {\rm grad}%
h\right\Vert ^{2}}{2h^{2}}g(X,X)\right) \\
&&+\frac{Xh}{h}g(X,X)\left( (n_{2}-1)\frac{\left\Vert gradf\right\Vert ^{2}}{%
f^{2}}-n_{2}\pi (P)-n_{2}\frac{1}{f}\triangle f\right) .
\end{eqnarray*}%
The proof is completed. $\blacksquare $

\begin{rem-new}
Einstein-like manifolds are natural extension of Einstein manifolds.
Einstein-like manifolds admitting different curvature conditions were
considered by Calvaruso \cite{Calvaruso}. Einstein-like manifolds of
dimension $3$ and $4$ are studied in \cite{Berndt, Bueken}. Projective
spaces and spheres furnished with class ${\cal A}$ or class ${\cal B}$
Einstein-like metrics were classified in \cite{Peng}. An interesting study
in \cite{Mantica} shows \ that Einstein-like Generalized Robertson-Walker
spacetimes are perfect fluid space-times except one class of Gray's
decomposition.
\end{rem-new}

\noindent {\bf Acknowledgement:} We thank the reviewer for valuable
suggestions and corrections.

\noindent Department of Mathematics \& Statistics\newline
School of Mathematical \& Physical Sciences\newline
Dr. Harisingh Gour University\newline
Sagar-470003, Madhya Pradesh\newline
punam\_2101@yahoo.co.in \medskip

and

\noindent Universit\'{e} Alioune Diop de Bambey\newline
UFR SATIC, D\'{e}partement de Math\'{e}matiques\newline
ER-ANLG\newline
B.P. 30, Bambey, S\'{e}n\'{e}gal\newline
abdoulsalam.diallo@uadb.edu.sn \medskip

\end{document}